\documentclass[12pt,a4paper,oneside,titlepage]{article}
\usepackage{amsmath,amsfonts,amssymb,amsthm,amscd}
\newtheorem{Theorem}{Theorem}[section]

\newtheorem{Example}{Example}[section]
\newcommand{\R}{\mathbb R}

\newcommand{\Z}{\mathbb Z}
\newcommand{\N}{\mathbb N}

\title{On a Theorem of Legendre on Diophantine Approximation}
\author{Jaroslav Han\v cl and Tho Phuoc Nguyen}

\begin{document}
	
	\maketitle
	\date{}
	
	\begin{abstract}
		Legendre's theorem states that every irreducible fraction $\dfrac{p}{q}$ which satisfies the inequality $\left |\alpha-\dfrac{p}{q} \right | < \dfrac{1}{2q^2}$ is convergent to $\alpha$. Later Barbolosi and Jager improved this theorem. In this paper we refine these results.
	\end{abstract}

	\section{Introduction} 
	
	The theory of simple continued fractions plays the most important tools in mathematical analysis, probability theory, physics, approximation theory and other branches of natural sciences. During the last three centuries many famous mathematicians came with interesting results in Diophantine approximations and continued fractions. Among them let us mention for example Borel \cite{borel2}, Dirichlet \cite{Dirichlet}, Hurwitz \cite{hurwitz}, Legendre \cite{legendre} or Vahlen \cite{vahlen}.
	
	In 1830, Adrien - Marie Legendre \cite{legendre}, page 23, proved the sufficient condition for a fraction $\dfrac{p}{q}$ to be convergent of a real number $\alpha$. 
	\begin{Theorem} \label{legendre.theorem}  
		(Legendre) Let $p$ and $q$ be relatively prime integers with $q>0$ and such that
		\begin{equation}
			\Bigl\lvert \alpha -\frac pq\Bigr\rvert< \frac 1 {2 q^2}\,.
		\end{equation} 
		Then $\dfrac pq$ is a convergent to $\alpha$. 
	\end{Theorem}
	In Theorem \ref{hanthot1} we refine this theorem, replacing $2$ by  $2-(q-1)/q^2$. The proof is elementary in character but lengthy, as it  involvies a detailed case analysis. Theorems 
	\ref{hanthot2} and \ref{hanthot3} give other alternative to Theorem \ref{hanthot1}. 
	
	Legendre's theorem has some history. In 1965, Billingsley \cite{billingsley} made use Legendre's method in the ergodic theory. In 1988, Ito \cite{Ito1988} tried to extend Legendre's constant $1/2$ to the other kind of continued fractions. To do this he used a special algorithms how to find the best constant. Such kind of methods are used for example in the theory of dynamical systems and ergodic theory. For other result see \cite{Itoosaka}, \cite{Ito1985} and \cite{nakada}.
	
	There is a nice result of Koksma \cite{koksmatheorem} stated by the following theorem
	\begin{Theorem} (Koksma) \label{koksma}
		If $\dfrac{p}{q}$ is a rational, and $\alpha$ an irrational number and if $q|q\alpha-p|<\dfrac{2}{3}$, then $\dfrac{p}{q}$ is either a convergent or a first mediant of $\alpha$. The constant $\dfrac{2}{3}$ is best possible.
	\end{Theorem}
	
	Barbolosi and Jager \cite{barbolosi1994} refined Legendre's theorem by considering two special cases from which the constant $2/3$ appears, they prove the following theorem.
	\begin{Theorem} \label{barbolosi}
		(Barbolosi and Jager) Let $\dfrac{p}{q}=[b_0;b_1,\dots,b_n]$ where $b_n\geq2$, $(p,q)=1,q>0$ be the simple continued fraction expansion and let $\alpha$ be an irrational number. If $(-1)^nsgn(\alpha-\dfrac{p}{q})=1$ and $q|q\alpha-p|<\dfrac{2}{3}$ then $\dfrac{p}{q}$ is a convergent of $\alpha$.
	\end{Theorem}
	This is an interesting result because it shows the possibility of improving the constant $1/2$ in Legendre's theorem. However, in some cases, determining the length of the continued fraction representing $p/q$ is not easy. 
	
	We improve  Theorem \ref{koksma} and give some alternative to Theorem \ref{barbolosi} in Theorem \ref{hanthot6}, when we replace $2/3$ by $(3/2-((q-1)(q-2))/(2q^3))^{-1}$. The similar method was used in \cite{bh} and \cite{hancl1}-\cite{hancl3} but for completely different results.
	
	Interesting results we can find in \cite{barbolosi1988}, \cite{grace}, \cite{jager&kraaikamp} and \cite{koksma1951}-\cite{kraaikamp}. An excellent source of basic background is  Hardy and Wright \cite{hardy}. The books of  Wall \cite{wall}, Hensley \cite{hensley}, Khinchin \cite{khinchin}, Karpenkov \cite{karpenkov}, Thron \cite{jones},  Rosen \cite{rosen}, Borwein and Borwein \cite{borwein} and Feldman and Nesterenko \cite{feldman} are also very useful.

	\section{Notation} 
	
	Throughout the paper, $\Z^+$, $\N_0$ and $\R$ will denote the sets of positive, non-negative integers and real numbers respectively. Let $\alpha$ be a real number and suppose $n\in\N_0$\,.  Let $\alpha=[a_0;a_1,a_2,\dots ]$ be  its simple continued fraction expansion.  Also let $\dfrac {p_n}{q_n}=[a_0;a_1,a_2,\dots ,a_n]$ be its $n$-th convergent. 
	The following recurrence relations for convergents are known
	\begin{align*}
		p_{-2}&=0\,, &  p_{-1}&=1\,, &  p_0&=a_0\,,\\
		q_{-2}&=1\,, &    q_{-1}&=0\,,   &  q_0&=1\,,
	\end{align*}
	$$p_1=a_1a_0+1\,,\hspace{2cm}  p_{n+2}=a_{n+2}p_{n+1}+p_n\,,$$
	$$q_1=a_1\,,   \hspace{2cm}          q_{n+2}=a_{n+2}q_{n+1}+q_n\,.$$
	For a simple continued fraction expansion  we have that
	$$\alpha=[a_0;a_1,a_2,\dots ]=[a_0;a_1,a_2,\dots ,a_n,[a_{n+1};a_{n+2},a_{n+3},\dots ]].$$
	Therefore, we can write
	$$\alpha=[a_0;a_1,a_2,\dots,a_n,r_{n+1}]$$
	where $r_{n+1}=[a_{n+1};a_{n+2},a_{n+3},\dots ]$. From this we obtain that
	\begin{equation} \label{hantho.eq.notation1}
		\alpha=\dfrac{r_{n+1}p_n+p_{n-1}}{r_{n+1}q_n+q_{n-1}}.
	\end{equation}
	Taking a difference of two consecutive convergents we obtain that
	\begin{equation} \label{hantho.eq.notation2}
		q_{n+1}p_n-p_{n+1}q_n=(-1)^{n+1}\,.
	\end{equation} 
	For the finite simple continued fraction expansion, if we have \\
	$\alpha=[a_0;a_1,a_2,\dots ,a_k]$ with $k\geq 1$\,, then we suppose that $a_k\geq2$\,.  
	
	The sequence of mediants of real irrational number $\alpha$ is the sequence of irreducible fractions of the form $\dfrac{bp_n+p_{n-1}}{bq_n+q_{n-1}}$ with $ n\geq0, b=1, 2, \dots, a_{n+1}-1$, ordered in such a way that the denominators form an increasing sequence. The mediant with $b=1$ or $b=a_{n+1}-1$ is called nearest mediant.
	
	More details on the discussion in this section can be found in \cite{barbolosi1994} or \cite{schmidt}, pages $7$ to $10$\,.

	\section{New Results}
	Our main result is the following theorem which concerns with the Legendre's theorem.
	\begin{Theorem} \label{hanthot1}  
		Let $p,q$ be relatively prime integers with $q\geq1$ and
		\begin{equation}\label{hantho.eq.int1}
			\left | \alpha -\frac{p}q{} \right |\leq\frac{1}{\left ( 2-\dfrac{q-1}{q^2} \right )q^{2}}.
		\end{equation}
		Then $\dfrac{p}{q}$ is a convergent to $\alpha$ excluding the case when $\alpha=[a_0;2]$ and $\dfrac{p}{q}=[a_0+1]$. For this special case, we have equality in \eqref{hantho.eq.int1} and $\dfrac{p}{q}$ is a nearest mediant of $\alpha$.
	\end{Theorem}
	This theorem we can improve in the following way.
	\begin{Theorem} \label{hanthot2}
		Let $p,q$ be relatively prime integers with $q\geq1$ and
		\begin{equation}\label{hantho.eq.int2}
			\left | \alpha -\frac{p}q{} \right |\leq\frac{1}{\left ( 2-\dfrac{1}{q} \right )q^{2}}.
		\end{equation}
		Then $\dfrac{p}{q}$ is a convergent to $\alpha$ excluding the cases
		\begin{enumerate}
			\item $\alpha=[a_0]$ and $\dfrac{p}{q}=[a_0-1]$.
			\item $\alpha=[a_0]$ and $\dfrac{p}{q}=[a_0+1]$.
			\item $\alpha=[a_0;a_1,a_2,\dots], (a_1\geq2)$ and $\dfrac{p}{q}=[a_0+1]$.
			\item $\alpha=[a_0;3]$ and $\dfrac{p}{q}=[a_0;2]$.
			\item $\alpha=[a_0;a_1,2]$ and $\dfrac{p}{q}=[a_0;a_1+1]$.
		\end{enumerate}  
		For all special cases 1-5, we have that $\dfrac{p}{q}$ is the nearest mediant of $\alpha$. For cases 1, 2, 4 and 5 we have equality in \eqref{hantho.eq.int2}, case 3 satisfies sharp inequality in \eqref{hantho.eq.int2}.
	\end{Theorem}
	For every large $n$ we have the following result.
	\begin{Theorem} \label{hanthot3}
		Let $\alpha$ be a real number, $p,q$ be relatively prime integers with $q\geq1$, $n$ be a positive integer and $\dfrac{p_{n-1}}{q_{n-1}}$ be a convergent of $\alpha$ and also $\dfrac{p}{q}$. Assume that
		\begin{equation}\label{hantho.eq.int3}
			\left | \alpha -\frac{p}q{} \right |\leq\frac{1}{\left ( 2-\dfrac{q_{n-1}}{q} \right )q^{2}}.
		\end{equation}
		Then $\dfrac{p}{q}$ is a convergent to $\alpha$, excluding the cases
		\begin{enumerate}
			\item $\alpha=[a_0;2]$ and $\dfrac{p}{q}=[a_0+1]$.
			\item $\alpha=[a_0;3]$ and $\dfrac{p}{q}=[a_0;2]$.
			\item $\alpha=[a_0;a_1,\dots,a_{n-1},a_n,2]$ and $\dfrac{p}{q}=[a_0;a_1,\dots,a_{n-1},a_n+1]$.
		\end{enumerate}
		For all these cases $\dfrac{p}{q}$ is the nearest mediant of $\alpha$ and we have equality in \eqref{hantho.eq.int3}.
	\end{Theorem}
	The next theorem improves the result of Barbolosi and Jager.
	\begin{Theorem}\label{hanthot6}
		Let $p,q$ be relatively prime integers with $q\geq1$. If
		\begin{equation} \label{hantho.eq.int6}
			\left|\alpha-\dfrac{p}{q} \right|\leq\dfrac{1}{\left(1-\dfrac{1}{2q}\right)q^2}
		\end{equation}
		then $\dfrac{p}{q}$ is either a convergent or nearest mediant of $\alpha$, excluding the case
		\begin{enumerate}
			\item $\alpha=[a_0]$ and $\dfrac{p}{q}=[a_0-2]$.
			\item $\alpha=[a_0;a_1,a_2,\dots]$ and $\dfrac{p}{q}=[a_0+2]$.
			\item $\alpha=[a_0;6]$ and $\dfrac{p}{q}=[a_0;2]$.
			\item $\alpha=[a_0;5,a_2,\dots]$ and $\dfrac{p}{q}=[a_0;2]$.
			\item $\alpha=[a_0;a_1,4]$ and $\dfrac{p}{q}=[a_0;a_1,2]$.
			\item $\alpha=[a_0;5]$ and $\dfrac{p}{q}=[a_0;3]$.
			\item $\alpha=[a_0;4,a_2,\dots]$ and $\dfrac{p}{q}=[a_0;2]$.
		\end{enumerate}
	\end{Theorem}
	\begin{Example} \label{example1}
		Let $\alpha=\displaystyle\sum_{n=1}^{\infty}\dfrac{1}{2^{2^n-1}A^{2^n}}$ where $A \in \Z^+$. Set\\
		$$\dfrac{p}{q}=\dfrac{\left(2^{2^N-1}A^{2^N}\right)\displaystyle\sum_{n=1}^{N}\dfrac{1}{2^{2^n-1}A^{2^n}}}{2^{2^N-1}A^{2^N}}=\displaystyle\sum_{n=1}^{N}\dfrac{1}{2^{2^n-1}A^{2^n}} $$
		where $N \in \Z^+$. From Theorem \ref{hanthot1} we obtain that $\dfrac{p}{q}$ is a convergent of $\alpha$ which is not an immediate consequence of Legendre's theorem.
	\end{Example}
	
	\begin{Example} \label{example2}
		From Example \ref{example1} we obtain that $\displaystyle\sum_{n=1}^{N}\dfrac{1}{2^{2^n-1}}$ is a convergent of  $\displaystyle\sum_{n=1}^{\infty}\dfrac{1}{2^{2^n-1}}$. This is not an immediate consequence of Legendre's Theorem.
	\end{Example}
	\begin{Example} \label{example4}
		Let $\alpha=\displaystyle\sum_{n=1}^{\infty}\dfrac{1}{2^{2^n}A^{2^n}}$ where $A \in \Z^+$. Set\\
		$$\dfrac{p}{q}=\dfrac{\left(2^{2^N}A^{2^N}\right)\displaystyle\sum_{n=1}^{N}\dfrac{1}{2^{2^n}A^{2^n}}}{2^{2^N}A^{2^N}}=\displaystyle\sum_{n=1}^{N}\dfrac{1}{2^{2^n}A^{2^n}} $$
		where $N \in \Z^+$. From Theorem \ref{hanthot6} we obtain that $\dfrac{p}{q}$ is a convergent or nearest mediant of $\alpha$ which is not an immediate consequence of Barbolosi and Jager's theorem.
	\end{Example}

	\section{Proofs}

	\begin{proof}[Proof of Theorem \ref{hanthot1}]  
		
		Let $\alpha =[a_0;a_1,a_2,\dots,a_{n-1},a_n,\dots ]$ be a simple continued fractional expansion of the number $\alpha$. For any irreducible fraction $\dfrac{p}{q}$ which is not a convergent of $\alpha$ and $m,n \in \N_0$ we can write as
		\[\dfrac{p}{q}=\left [ a_0;a_1,a_2,\dots,a_{n-1},b_n,b_{n+1},\dots,b_{n+m} \right ]\]
		where $b_n,b_{n+1},\dots,b_{n+m} \in \Z^+$ and $b_{n+m}\geq2$.
		\begin{enumerate}
			\item 
			Suppose that $q=1$. Then $\dfrac{p}{q}=[b_0]=\dfrac{b_0}{1}$. Now we prove that
			$$\left|\alpha-\dfrac{p}{q} \right|>\dfrac{1}{\left(2-\dfrac{q-1}{q^2}\right)q^2}=\dfrac{1}{2}.$$
			\begin{enumerate}
				\item Assume that $a_0=b_0$ then $\dfrac{p}{q}$ is a convergent of $\alpha$.
				\item Let $a_0>b_0$. Then we obtain that
				$$\left|\alpha-\dfrac{p}{q} \right|=a_0-b_0+[0;a_1,\dots] \geq 1+[0;a_1,\dots]>\dfrac{1}{2}. $$
				\item Assume that $b_0 \geq a_0+2$. Then we have
				$$\left|\alpha-\dfrac{p}{q} \right|=b_0-a_0-[0;a_1,\dots] \geq 2-[0;a_1,\dots]>\dfrac{1}{2}.$$
				\item Suppose that $b_0=a_0+1$. It yields
				$$\left|\alpha-\dfrac{p}{q} \right|=b_0-a_0-[0;a_1,\dots] = 1-[0;a_1,\dots].$$
				\begin{enumerate}
					\item If $a_1$ does not exist then we obtain that
					$$\left|\alpha-\dfrac{p}{q} \right|=b_0-a_0 = 1>\dfrac{1}{2}.$$
					\item If $a_1=1$ and $a_2$ exists then $\dfrac{p}{q}$ is a convergent of $\alpha$.
					\item If $a_1 \geq 2$ then we have
					$$\left|\alpha-\dfrac{p}{q} \right|=1-\dfrac{1}{a_1+[0;a_2,\dots]} \geq \dfrac{1}{2}.$$
					The equality occurs only in the case when $a_1=2$ and $a_2$ does not exist. This is the exception mentioned in Theorem \ref{hanthot1}.
				\end{enumerate}
			\end{enumerate}
			\item
			Now we suppose that $q\geq2$. Then we have $2-\dfrac{q-1}{q^2}>2-\dfrac{1}{q}$ for all $q\in \Z^+$. Therefore instead of \eqref{hantho.eq.int1} it is enough to prove \eqref{hantho.eq.int2}. The proof falls into two main cases. Here is the plan of our proof.\\
			a. $n=0$. \hspace{1cm}i. $a_0>b_0$. \\
			\phantom{1. $n=0$. \hspace{1cm}}ii. $b_0>a_0$. \\
			b. $n\geq1$. \hspace{1cm}i. $a_n>b_n$. \hspace{1cm}A. $d\geq1$ and $c\geq2$.\\
			\phantom{2. $n\geq1$. \hspace{1cm}a. $a_n>b_n$. \hspace{1cm}}B. $d=0$ and $c=1$.\\
			\phantom{2. $n\geq1$. \hspace{1cm}}ii. $b_n>a_n$. \hspace{1cm}A. $d\geq1$ and $c\geq2$.\\
			\phantom{2. $n\geq1$. \hspace{1cm}a. $a_n>b_n$. \hspace{1cm}}B. $d=0$ and $c=1$.
			\begin{enumerate}

				%
				%
				\item 
				Assume that $n=0$. Then we have
				$$\dfrac{p}{q}=[b_0;b_1,b_2,\dots,b_{m}]=\dfrac{cb_0+d}{c}$$
				where $\dfrac{d}{c}=[0;b_1,b_2,\dots,b_m]$ and $c>d>0$. Note that $c \neq 1$ otherwise $c=q=1$. Set $A=[0;a_1,\dots]$ then we have
				\begin{equation} \label{hantho.eq.case1}
					\left|\alpha-\dfrac{p}{q} \right|=\left|a_0-b_0+A-\dfrac{d}{c}\right|=\dfrac{1}{c^2}c^2\left|a_0-b_0+A-\dfrac{d}{c}\right|.
				\end{equation}
				Now we prove that
				$$\left|\alpha-\dfrac{p}{q}\right|>\dfrac{1}{\left(2-\dfrac{1}{q}\right)q^2}.$$
				From this, the fact that $q=c\geq2$ and \eqref{hantho.eq.case1} we obtain that it is enough to prove that
				$$\dfrac{1}{c^2}c^2\left|a_0-b_0+A-\dfrac{d}{c}\right|>\dfrac{1}{\left(2-\dfrac{1}{c}\right)c^2}$$
				which is equivalent to
				\begin{equation} \label{hantho.eq.new.case3}
					c(2c-1)\left|a_0-b_0+A-\dfrac{d}{c}\right|>1.
				\end{equation}
				Let us suppose two cases.
				\begin{enumerate}
					\item 
					Assume that $a_0>b_0$. This yields that inequality \eqref{hantho.eq.new.case3} has the form
					$$c(2c-1)\left(a_0-b_0+A-\dfrac{d}{c}\right)>1$$
					which is obviously true since $c\geq2$ and\\
					\vspace{0.2cm}$c(2c-1)\left(a_0-b_0+A-\dfrac{d}{c}\right)\geq c(2c-1)\left(1+0-\dfrac{d}{c}\right)$\\
					\vspace{0.2cm}\phantom{$c(2c-1)\left(a_0-b_0+A-\dfrac{d}{c}\right)$} $=(2c-1)\left(c-d\right)\geq3>1$.\\
					Hence inequality \eqref{hantho.eq.new.case3} follows. 
					\item
					Let $b_0>a_0$. Then inequality \eqref{hantho.eq.new.case3} is equivalent to
					$$c(2c-1)\left(b_0-a_0-A+\dfrac{d}{c}\right)>1$$
					which is also obviously true since $c\geq2,d\geq1$ and\\
					\vspace{0.2cm}$c(2c-1)\left(b_0-a_0-A+\dfrac{d}{c}\right)> c(2c-1)\left(1-1+\dfrac{d}{c}\right)$\\
					\vspace{0.2cm}\phantom{$c(2c-1)\left(a_0-b_0+A-\dfrac{d}{c}\right)$} $=(2c-1)d>1$.\\
					Therefore inequality \eqref{hantho.eq.new.case3} follows.

				\end{enumerate}
				%
				%
				%
				\item 
				Suppose that $n\geq1$. Then we have\\
				$\alpha=[a_0;a_1,a_2,\dots,a_{n-1},a_n,a_{n+1},\dots]
				=\left[a_0;a_1,\dots,a_{n-1},a_n+r\right]$\\
				where $r=[0;a_{n+1},a_{n+2},\dots]$ satisfies $0\leq r<1$. If $r=0$ then $a_{n+1}$ doesn't exist. This and equality \eqref{hantho.eq.notation1} yield
				\begin{equation} \label{hantho.eq.case2}
					\alpha=\dfrac{(a_n+r)p_{n-1}+p_{n-2}}{(a_n+r)q_{n-1}+q_{n-2}}
				\end{equation}
				and
				$$\dfrac{p}{q}=[a_0;a_1,a_2,\dots,a_{n-1},b_n,b_{n+1},\dots,b_{n+m}].$$
				Set $\dfrac{d}{c}=[0;b_{n+1}\dots,b_{n+m}]$ where $c>d$. If $d=0$, then $c=1$. Otherwise $p$ and $q$ are not coprime. This and equality \eqref{hantho.eq.notation1} imply 
				$$\dfrac{p}{q}=\dfrac{p_{n-1}\left(b_n+\dfrac{d}{c}\right)+p_{n-2}}{q_{n-1}\left(b_n+\dfrac{d}{c}\right)+q_{n-2}}=\dfrac{p_{n-1}(cb_n+d)+cp_{n-2}}{q_{n-1}(cb_n+d)+cq_{n-2}}.$$
				From this, \eqref{hantho.eq.notation2} and \eqref{hantho.eq.case2} we obtain that\\
				\vspace{0.2cm}$\left|\alpha-\dfrac{p}{q} \right|=\left| \dfrac{(a_n+r)p_{n-1}+p_{n-2}}{(a_n+r)q_{n-1}+q_{n-2}} - \dfrac{p_{n-1}(cb_n+d)+cp_{n-2}}{q_{n-1}(cb_n+d)+cq_{n-2}} \right|$ \\
				\vspace{0.2cm}\phantom{$\left|\alpha-\dfrac{p}{q} \right|$} $=\dfrac{|cb_n+d-ca_n-cr|}{((a_n+r)q_{n-1}+q_{n-2})(q_{n-1}(cb_n+d)+cq_{n-2})}$ \\
				\vspace{0.2cm}\phantom{$\left|\alpha-\dfrac{p}{q} \right|$} $=\dfrac{1}{q^2}\dfrac{q^2|cb_n+d-ca_n-cr|}{((a_n+r)q_{n-1}+q_{n-2})(q_{n-1}(cb_n+d)+cq_{n-2})}$ \\
				\vspace{0.2cm}\phantom{$\left|\alpha-\dfrac{p}{q} \right|$} $=\dfrac{1}{q^2}\dfrac{q|cb_n+d-ca_n-cr|}{(a_n+r)q_{n-1}+q_{n-2}}$ \\
				\vspace{0.2cm}\phantom{$\left|\alpha-\dfrac{p}{q} \right|$} $=\dfrac{1}{q^2}\dfrac{1}{\dfrac{(a_n+r)q_{n-1}+q_{n-2}}{q|cb_n+d-ca_n-cr|}}$ \\
				\vspace{0.2cm}\phantom{$\left|\alpha-\dfrac{p}{q} \right|$} $=\dfrac{1}{q^2}\dfrac{c}{\dfrac{c(a_n+r)q_{n-1}+cq_{n-2}-(cb_n+d)q_{n-1}+(cb_n+d)q_{n-1}}{q|cb_n+d-ca_n-cr|}}$ \\
				\vspace{0.2cm}\phantom{$\left|\alpha-\dfrac{p}{q} \right|$} $=\dfrac{1}{q^2}\dfrac{c}{\dfrac{q_{n-1}(c(a_n+r)-cb_n-d)+q}{q|cb_n+d-c(a_n+r)|}}$
				\begin{equation} \label{hantho.eq.new.case1method2}
					=\dfrac{1}{q^2}\dfrac{c}{\dfrac{1}{|c(b_n-a_n-r)+d|}+\dfrac{q_{n-1}sgn(a_n-b_n)}{q}}.
				\end{equation}
				\begin{enumerate}
					\item 
					If $a_n>b_n$ then \eqref{hantho.eq.new.case1method2} has the form
					\begin{equation} \label{hantho.eq.case.2a}
						\left|\alpha-\dfrac{p}{q} \right|=\dfrac{1}{q^2}\dfrac{c}{\dfrac{1}{c(a_n-b_n+r)-d}+\dfrac{q_{n-1}}{q}}.
					\end{equation}
					\begin{enumerate}
						\item 
						Assume that $d\geq1$. Then $c\geq2$. From this, \eqref{hantho.eq.case.2a} and the fact that $q=q_{n-1}(cb_n+d)+cq_{n-2}$ we obtain that\\
						\vspace{0.2cm}$\left|\alpha-\dfrac{p}{q} \right|=\dfrac{1}{q^2}\dfrac{c}{\dfrac{1}{c(a_n-b_n+r)-d}+\dfrac{q_{n-1}}{q}}\geq \dfrac{1}{q^2}\dfrac{2}{1+\dfrac{q_{n-1}}{q}}$\\
						\vspace{0.2cm}\phantom{$\left|\alpha-\dfrac{p}{q} \right|$}
						$= \dfrac{1}{q^2}\dfrac{2}{1+\dfrac{q_{n-1}}{q_{n-1}(cb_n+d)+cq_{n-2}}}$\\
						\begin{equation} \label{hantho.eq.t1.case2ai}
							\hspace{0.9cm}= \dfrac{1}{q^2}\dfrac{2}{1+\dfrac{1}{cb_n+d+c\dfrac{q_{n-2}}{q_{n-1}}}}\geq\dfrac{1}{q^2}\dfrac{3}{2}>\dfrac{1}{q^2}\dfrac{1}{2-\dfrac{1}{q}}.
						\end{equation}
						%
						\item 
						Suppose that $d=0$. Then $c=1$ and $b_{n+1}$ doesn't exist. Therefore $q=b_nq_{n-1}+q_{n-2}$ where $b_n\geq2$. From this and \eqref{hantho.eq.case.2a} we obtain that
						\begin{align} \label{hantho.eq.t1.case2aii}
							\left|\alpha-\dfrac{p}{q} \right|&=\dfrac{1}{q^2}\dfrac{c}{\dfrac{1}{c(a_n-b_n+r)-d}+\dfrac{q_{n-1}}{q}} \nonumber \\
							&\geq \dfrac{1}{q^2}\dfrac{1}{1+\dfrac{q_{n-1}}{q}}.
						\end{align}
						Hence\\
						\vspace{0.2cm} 
						$\left|\alpha-\dfrac{p}{q} \right| \geq \dfrac{1}{q^2}\dfrac{1}{1+\dfrac{q_{n-1}}{q}} =\dfrac{1}{q^2}\dfrac{1}{1+\dfrac{q_{n-1}}{b_nq_{n-1}+q_{n-2}}} $\\
						\vspace{0.2cm} \phantom{$\left|\alpha-\dfrac{p}{q} \right|$}
						$=\dfrac{1}{q^2}\dfrac{1}{1+\dfrac{1}{b_n+\dfrac{q_{n-2}}{q_{n-1}}}}\geq \dfrac{1}{q^2}\dfrac{1}{1+\dfrac{1}{b_n}} \geq \dfrac{1}{q^2}\dfrac{2}{3} \geq \dfrac{1}{q^2}\dfrac{1}{2-\dfrac{1}{q}}.$
					\end{enumerate}
					%
					\item 
					Let $b_n>a_n$. Then \eqref{hantho.eq.new.case1method2} has the form
					\begin{equation} \label{hantho.eq.case.2b}
						\left|\alpha-\dfrac{p}{q} \right|=\dfrac{1}{q^2}\dfrac{c}{\dfrac{1}{c(b_n-a_n-r)+d}-\dfrac{q_{n-1}}{q}}.
					\end{equation}
					\begin{enumerate}
						\item 
						Assume that $d\geq1$. Then $c\geq2$. From this, \eqref{hantho.eq.case.2b} and the facts $0\leq r<1$, $q_{n-1}\geq1$ we obtain that
						\begin{align} \label{hantho.eq.t1.case2bi}
							\left|\alpha-\dfrac{p}{q} \right|&=\dfrac{1}{q^2}\dfrac{c}{\dfrac{1}{c(b_n-a_n-r)+d}-\dfrac{q_{n-1}}{q}}\nonumber \\ 
							&> \dfrac{1}{q^2}\dfrac{2}{\dfrac{1}{2(1-1)+1}-\dfrac{q_{n-1}}{q}} \nonumber\\
							&= \dfrac{1}{q^2}\dfrac{2}{1-\dfrac{q_{n-1}}{q}}\geq \dfrac{1}{q^2}\dfrac{2}{1-\dfrac{1}{q}}> \dfrac{1}{q^2}\dfrac{1}{2-\dfrac{1}{q}}.
						\end{align}
						%
						%
						\item 
						Suppose that $d=0$. Then $c=1$ and $b_{n+1}$ doesn't exist. Therefore $q=b_nq_{n-1}+q_{n-2}$ where $b_n\geq2$. Now we prove that
						\begin{equation} \label{hantho.eq.case.2b.sp}
							b_n-a_n-r\geq\dfrac{1}{2}.
						\end{equation}
						Let $b_n=a_n+1$.\\
						$\bullet$ If $a_{n+1}$ does not exist then $b_n-a_n-r=1>\dfrac{1}{2}$.\\
						$\bullet$ If $a_{n+1}=1$ then $\dfrac{p}{q}$ is a convergent of $\alpha$.\\
						$\bullet$ If $a_{n+1}=2$ and $a_{n+2}$ does not exist then $r=\dfrac{1}{2}$ and we have
						$$b_n-a_n-r=1-\dfrac{1}{2}=\dfrac{1}{2}.$$
						$\bullet$ If $a_{n+1}=2$ and $a_{n+2}\geq1$ then
						$$b_n-a_n-r = 1-[0;2,a_{n+2},\dots]>\dfrac{1}{2}.$$
						$\bullet$ If $a_{n+1} \geq 3$ then we have $b_n-a_n-r \geq 1-[0;3]=\dfrac{2}{3}>\dfrac{1}{2}.$\\
						Let $b_n \neq a_n+1$ then $b_n-a_n-r \geq 2-r > \dfrac{1}{2}$. Hence \eqref{hantho.eq.case.2b.sp} follows.\\
						From \eqref{hantho.eq.case.2b}, \eqref{hantho.eq.case.2b.sp}, the facts that $c=1$, $d=0$ and $q_{n-1}\geq1$ we obtain that\\
						\vspace{0.2cm}$\left|\alpha-\dfrac{p}{q} \right|=\dfrac{1}{q^2}\dfrac{c}{\dfrac{1}{c(b_n-a_n-r)+d}-\dfrac{q_{n-1}}{q}} $\\
						\phantom{$\left|\alpha-\dfrac{p}{q} \right|$} $\geq \dfrac{1}{q^2}\dfrac{1}{2-\dfrac{q_{n-1}}{q}}\geq\dfrac{1}{q^2}\dfrac{1}{2-\dfrac{1}{q}}$.
					\end{enumerate}
				\end{enumerate}
			\end{enumerate}
		\end{enumerate}
		The proof of Theorem \ref{hanthot1} is complete.
	\end{proof}
	
	\begin{proof}[Proof of Theorem \ref{hanthot2}]
		We only follow Case 2 of Theorem \ref{hanthot1} with following discussions.
		\begin{enumerate}
			\item 
			In case 2(a)i we proved that $\left|\alpha-\dfrac{p}{q} \right|>\dfrac{1}{\left(2-\dfrac{1}{q}\right)q^2}$ for all cases when $q=c\geq2$. Suppose that $q=c=1$ then $d=0$ and we have
			$$c(2c-1)\left(a_0-b_0+A-\dfrac{d}{c}\right)= a_0-b_0+A\geq1.$$
			The equality occurs when $b_0=a_0-1$ and $a_1$ does not exist, so $A=0$. This is the first exception.
			\item  
			Suppose that $q=c=1,d=0$ in case 2(a)ii then we have
			$$c(2c-1)\left(b_0-a_0-A+\dfrac{d}{c}\right)=b_0-a_0-A.$$
			\begin{enumerate}
				\item Let $b_0\geq a_0+2$. Then we obtain that
				$$b_0-a_0-A \geq 2-A>1.$$
				\item Assume that $b_0= a_0+1$. It implies that
				$$b_0-a_0-A=1-A.$$
				$\bullet$ If $a_1$ does not exist then $A=0$ and we have $b_0-a_0-A=1-A=1$. This is the second exception.\\
				$\bullet$ Let $a_1=1$ and $a_2$ exists. So $\alpha=[a_0;1;a_2,\dots]$ and $\dfrac{p}{q}=[a_0+1]$ is convergent of $\alpha$.\\
				$\bullet$ If $a_1 \geq 2$ then we obtain that
				$$b_0-a_0-A=1-\dfrac{1}{a_1+[0;a_2,\dots]}<1$$ 
				which is the third exception.
			\end{enumerate}
			\item In case 2(b)iB we proved that $\left|\alpha-\dfrac{p}{q} \right|\geq\dfrac{1}{\left(2-\dfrac{1}{q}\right)q^2}$.
			The equality occurs only when $n=1,b_1=2,a_1=3$ and $a_2$ does not exist. Hence $r=0$. This is the fourth exception.
			\item In case 2(b)iiB, we also proved that $\left|\alpha-\dfrac{p}{q} \right|\geq\dfrac{1}{\left(2-\dfrac{1}{q}\right)q^2}$. The equality occurs when $n=1,b_1=a_1+1,a_{n+1}=2$ and $a_{n+2}$ does not exist. This is the fifth exception.
		\end{enumerate}
		Proofs of other cases are the same like in the proof of Theorem \ref{hanthot1}.\\
		The proof of Theorem \ref{hanthot2} is complete.
	\end{proof}
	
	\begin{proof}[Proof of Theorem \ref{hanthot3}]
		The proof of this theorem follows case 2 in the proof of Theorem \ref{hanthot1} with some following discussions.
		\begin{enumerate}
			\item 
			For case 2(a) if $n=0$ then $q_{n-1}=q_{-1}=0$. So it is enough to prove that
			$$\left|\alpha-\dfrac{p}{q} \right|=\dfrac{1}{c^2}c^2\left|a_0-b_0+A-\dfrac{d}{c} \right|>\dfrac{1}{\left(2-\dfrac{q_{n-1}}{q}\right)q^2}=\dfrac{1}{2q^2}=\dfrac{1}{2c^2}$$
			which is equivalent to
			\begin{equation} \label{hantho.eq.t3.1}
				2c^2\left|a_0-b_0+A-\dfrac{d}{c} \right|>1.
			\end{equation}
			\begin{enumerate}
				\item Suppose that $a_0>b_0$ then inequality \eqref{hantho.eq.t3.1} has the form
				$$2c^2\left(a_0-b_0+A-\dfrac{d}{c} \right)>1$$
				which is obviously true since
				$$2c^2\left(a_0-b_0+A-\dfrac{d}{c} \right)\geq2c(c-d)\geq2.$$
				This is also true for $q=c\geq1$ and $d\geq0$.
				\item Assume that $b_0\geq a_0+2$ then inequality \eqref{hantho.eq.t3.1} is equivalent to
				$$2c^2\left(b_0-a_0-A+\dfrac{d}{c} \right)>1$$
				which is obviously true since
				$$2c^2\left(b_0-a_0-A+\dfrac{d}{c} \right)>2c^2\left(2-1+\dfrac{d}{c} \right) \geq 2>1$$
				for all values of $c\geq1,d \geq 0$.
				\item Let $b_0=a_0+1$ then inequality \eqref{hantho.eq.t3.1} is equivalent to
				$$2c^2\left(b_0-a_0-A+\dfrac{d}{c} \right)>1$$
				which is obviously true since
				$$2c^2\left(b_0-a_0-A+\dfrac{d}{c} \right)>2c^2\left(1-1+\dfrac{d}{c} \right) \geq 4>1$$
				for all values of $c\geq2, d\geq1$.\\
				Now we suppose that $c=1,d=0$ then we obtain 
				$$2c^2\left(b_0-a_0-A+\dfrac{d}{c} \right)=2(1-A).$$
				$\bullet$ If $a_1$ does not exist then $A=0$ and we have $2(1-A)=2>1.$\\
				$\bullet$ If $a_1=1$ and $a_2$ exists then $\dfrac{p}{q}$ is a convergent of $\alpha$.\\
				$\bullet$ If $a_1\geq2$ then we have $2(1-A) \geq 2\left(1-\dfrac{1}{a_1}\right)\geq1.$ The equality occurs when $\alpha=[a_0;2]$ and $\dfrac{p}{q}=[a_0+1]$, which is the first exception.
			\end{enumerate}
			\item 
			Suppose that $n\geq1$ then we follow case 2(b) in the proof of Theorem \ref{hanthot1} with these exceptions.
			\begin{enumerate}
				\item 
				In case 2(b)iA we have 
				$$\left|\alpha-\dfrac{p}{q} \right| \geq \dfrac{1}{q^2}\dfrac{2}{1+\dfrac{q_{n-1}}{q}}>\dfrac{1}{q^2}\dfrac{1}{2-\dfrac{q_{n-1}}{q}}.$$
				\item 
				In case 2(b)iB we have $b_n\geq2$, $q=b_nq_{n-1}+q_{n-2}\geq2q_{n-1}$ and
				$$\left|\alpha-\dfrac{p}{q} \right| \geq \dfrac{1}{q^2}\dfrac{1}{1+\dfrac{q_{n-1}}{q}}\geq\dfrac{1}{q^2}\dfrac{1}{2-\dfrac{q_{n-1}}{q}}.$$
				The equality occurs when $\alpha=[a_0;3]$ and $\dfrac{p}{q}=[a_0;2]$, this is the second exception.
				\item 
				In case 2(b)iiA we have
				$$\left|\alpha-\dfrac{p}{q} \right| \geq \dfrac{1}{q^2}\dfrac{2}{1-\dfrac{q_{n-1}}{q}}>\dfrac{1}{q^2}\dfrac{1}{2-\dfrac{q_{n-1}}{q}}.$$
				\item 
				In case 2(b)iiB we have $\left|\alpha-\dfrac{p}{q} \right| \geq\dfrac{1}{q^2}\dfrac{1}{2-\dfrac{q_{n-1}}{q}}.$\\
				The equality occurs when $\alpha=[a_0;a_1,\dots,a_{n-1},a_n,2]$ \\
				and $\dfrac{p}{q}=[a_0;a_1,\dots,a_{n-1},a_n+1]$, this is the third exception.
			\end{enumerate}
		\end{enumerate}
		The proof of Theorem \ref{hanthot3} is complete.
	\end{proof}

	\begin{proof}[Proof of Theorem \ref{hanthot6}]
		Let $\alpha=[a_0;a_1,a_2,\dots,a_{n-1},a_n,a_{n+1},\dots]$ be a simple continued fraction expansion of number $\alpha$. For any irreducible fraction $\dfrac{p}{q}$ which is neither convergent nor nearest mediant of $\alpha$ and $m,n\in \N_0$ we can write as
		$$\dfrac{p}{q}=[a_0;a_1,a_2,\dots,a_{n-1},b_n,b_{n+1},\dots,b_{n+m}]$$
		where $b_n,b_{n+1},\dots,b_{n+m} \in \Z^+$, $b_n\neq a_n$ and $b_{n+m}\geq2$. If $b_{n+1}$ does not exist then $b_n \notin \{a_n+1,a_n-1\}$.
		
		The proof falls into two main cases. Here is the plan of our proof.\\
		1. $n=0.$ \hspace{1cm} a. $c\geq2,d\geq1$ and $a_0\geq b_0+1$. \\
		\phantom{1. $n=0.$ \hspace{1cm} }b. $ c=1,d=0$ and $a_0 \geq b_0+2$.\\
		\phantom{1. $n=0.$ \hspace{1cm} }c. $c\geq2,d\geq1$ and $b_0\geq a_0+1$.\\
		\phantom{1. $n=0.$ \hspace{1cm} }b. $ c=1,d=0$ and $b_0 \geq a_0+2$.\\
		2. $n\geq1.$ \hspace{1cm} a. $a_n>b_n.$ \hspace{1.7cm} i. $c\geq2,d\geq1$ and $a_n \geq b_n+1$.\\
		\phantom{2. $n\geq1.$ \hspace{1cm} a. $a_n>b_n.$ \hspace{1.7cm} }ii. $a_n>b_n,c=1,d=0.$\\
		\phantom{2. $n\geq1.$ \hspace{1cm} a. $a_n>b_n.$ \hspace{1.7cm} ii. $c=1$}A. $a_n \geq b_n+4.$\\
		\phantom{2. $n\geq1.$ \hspace{1cm} a. $a_n>b_n.$ \hspace{1.7cm} ii. $c=1$}B. $a_n=b_n+3,n\geq2.$\\
		\phantom{2. $n\geq1.$ \hspace{1cm} a. $a_n>b_n.$ \hspace{1.7cm} ii. $c=1$}C. $a_n=b_n+3,n=1.$\\
		\phantom{2. $n\geq1.$ \hspace{1cm} a. $a_n>b_n.$ \hspace{1.7cm} ii. $c=1$}D. $a_n=b_n+2,n\geq2.$\\
		\phantom{2. $n\geq1.$ \hspace{1cm} a. $a_n>b_n.$ \hspace{1.7cm} ii. $c=1$}E. $a_n=b_n+2,n=1.$\\
		\phantom{2. $n\geq1.$ \hspace{1cm} }b. $b_n \geq a_n+2, c\geq1$ and $d\geq0.$\\
		\phantom{2. $n\geq1.$ \hspace{1cm} a. $a_n>b_n.$ \hspace{1.7cm} }i. $c\geq2,d\geq1$ and $b_n \geq a_n+1.$\\
		\phantom{2. $n\geq1.$ \hspace{1cm} a. $a_n>b_n.$ \hspace{1.7cm} }ii. $c=1,d=0$ and $b_n \geq a_n+2.$
		\begin{enumerate}
			\item Suppose that $n=0$. Then $$\alpha=[a_0;a_1,a_2,\dots,a_{n-1},a_n,a_{n+1},\dots]=[a_0+A]$$
			and
			$$\dfrac{p}{q}=[b_0;b_1,\dots,b_m]=b_0+\dfrac{d}{c}=\dfrac{cb_0+d}{c}$$
			where $ A=[0;a_1,a_2,\dots]\in[0;1)$ and $\dfrac{d}{c}=[0;b_1,b_2,\dots,b_m]\in[0;1)$. Note that $b_0 \neq a_0,b_m\geq2$ and if $d=0,c=1$ which mean $b_{1}$ does not exist then $b_0 \notin \{a_0+1,a_0-1\}$.\\
			From this we have
			$$\left|\alpha-\dfrac{p}{q} \right|=\left|a_0-b_0+A-\dfrac{d}{c}\right|=\dfrac{1}{c^2}c^2\left|a_0-b_0+A-\dfrac{d}{c}\right|.$$
			Now we prove that
			$$\left|\alpha-\dfrac{p}{q} \right|>\dfrac{1}{\left(1-\dfrac{1}{2q}\right)q^2}$$
			which is equivalent to
			\begin{equation} \label{eq.t6.n=0}
				\left(c^2-\dfrac{c}{2}\right)\left|a_0-b_0+A-\dfrac{d}{c} \right|>1.
			\end{equation}
			\begin{enumerate}
				\item Let $c\geq2,d\geq1$ and $a_0\geq b_0+1$. It yields that inequality \eqref{eq.t6.n=0} can be written as
				$$\left(c^2-\dfrac{c}{2}\right)\left(a_0-b_0+A-\dfrac{d}{c} \right)>1$$
				which is obviously true since
				\begin{align*}
					\left(c^2-\dfrac{c}{2}\right)\left(a_0-b_0+A-\dfrac{d}{c} \right)&\geq\left(c^2-\dfrac{c}{2}\right)\left(1-\dfrac{d}{c} \right)\\
					&=\left(c-\dfrac{1}{2}\right)\left(c-d \right)>1.
				\end{align*}
				
				\item Assume that $c=1,d=0$ and $a_0\geq b_0+2$. It yields that inequality \eqref{eq.t6.n=0} can be written as
				$$\left(c^2-\dfrac{c}{2}\right)\left(a_0-b_0+A-\dfrac{d}{c} \right)>1$$
				which is obviously true since
				\begin{align*}
					\left(c^2-\dfrac{c}{2}\right)\left(a_0-b_0+A-\dfrac{d}{c} \right)&\geq \dfrac{1}{2}(a_0-b_0+A)\\
					&\geq \dfrac{1}{2}(2+A) \geq1.
				\end{align*}
				The equality occurs when $A=0$ and $a_0= b_0+2$. It implies that $\alpha=[a_0]$ and $\dfrac{p}{q}=[a_0-2]$. This is the first exception.
				
				\item Suppose that $c \geq 2,d \geq 1$ and $b_0 \geq a_0+1$. Then inequality \eqref{eq.t6.n=0} has the form
				$$\left(c^2-\dfrac{c}{2}\right)\left(b_0-a_0-A+\dfrac{d}{c} \right)>1$$
				which is obviously true since
				\begin{align*}
					\left(c^2-\dfrac{c}{2}\right)\left(b_0-a_0-A+\dfrac{d}{c} \right)&>\left(c^2-\dfrac{c}{2}\right)\left(1-1+\dfrac{d}{c} \right)\\
					&=\left(c^2-\dfrac{c}{2}\right)\dfrac{d}{c} =\left(c-\dfrac{1}{2}\right)d>1.
				\end{align*}
				
				\item Assume that $b_0 \geq a_0+2,c=1$ and $d=0$. It implies that inequality \eqref{eq.t6.n=0} is equivalent to
				$$\dfrac{1}{2}(b_0-a_0-A)>1.$$
				which is obviously true since\\
				$\bullet$ If $b_0 \geq a_0+3$ then we obtain that
				$$\dfrac{1}{2}(b_0-a_0-A) \geq \dfrac{1}{2}(3-A)>1.$$
				$\bullet$ Suppose that $b_0 = a_0+2$. Hence
				$$\dfrac{1}{2}(b_0-a_0-A)=\dfrac{1}{2}(2-A)\leq1.$$
				The equality occurs when $A=0$ which mean $a_1$ does not exist. On the other side, when $a_1$ exists the sharp inequality satisfied. These are the cases included in exception 2.
			\end{enumerate}
			\item Let $n\geq1$. Then $$\alpha=[a_0;a_1,\dots,a_{n-1},a_n,a_{n+1},\dots]=\dfrac{(a_n+r)p_{n-1}+p_{n-2}}{(a_n+r)q_{n-1}+q_{n-2}}$$ and
			$$\dfrac{p}{q}=[a_0;a_1,\dots,a_{n-1},b_n,b_{n+1},\dots,b_{n+m}]=\dfrac{(cb_n+d)p_{n-1}+cp_{n-2}}{(cb_n+d)q_{n-1}+cq_{n-2}}$$
			where $ r=[0;a_{n+1},\dots]\in [0;1)$ and
			$\dfrac{d}{c}=[0;b_{n+1},\dots,b_{n+m}]\in [0;1)$, $c>d,b_n \neq a_n$. If $d=0$ then $c=1$ otherwise $p,q$ are not coprime and then $b_n \notin \{a_n+1,a_n-1\}$. From this and \eqref{hantho.eq.new.case1method2} we obtain that
			$$\left|\alpha-\dfrac{p}{q} \right|=\dfrac{1}{q^2}\dfrac{c}{\dfrac{1}{|c(b_n-a_n-r)+d|}+\dfrac{q_{n-1}sgn(a_n-b_n)}{q}}.$$
			Now we prove that
			$$\left|\alpha-\dfrac{p}{q} \right|>\dfrac{1}{\left(1-\dfrac{1}{2q}\right)q^2}$$
			which is equivalent to
			\begin{equation} \label{eq.t6.case2}
				\dfrac{c}{\dfrac{1}{|c(b_n-a_n-r)+d|}+\dfrac{q_{n-1}sgn(a_n-b_n)}{q}}>\dfrac{1}{1-\dfrac{1}{2q}}.
			\end{equation}
			\begin{enumerate}
				\item Assume that $a_n>b_n$ then \eqref{eq.t6.case2} has the form
				\begin{equation} \label{eq.t6.case2a}
					\dfrac{c}{\dfrac{1}{c(a_n-b_n+r)-d}+\dfrac{q_{n-1}}{q}}>\dfrac{1}{1-\dfrac{1}{2q}}.
				\end{equation}
				Now we consider some cases
				\begin{enumerate}
					\item Let $a_n\geq b_n+1,c\geq 2$ and $d\geq1$. \\
					Then $q=(cb_n+d)q_{n-1}+cq_{n-2}\geq cb_n+d \geq3$. From this we obtain that
					\begin{align*}
						\dfrac{c}{\dfrac{1}{c(a_n-b_n+r)-d}+\dfrac{q_{n-1}}{q}}&\geq \dfrac{c}{\dfrac{1}{c-d}+\dfrac{q_{n-1}}{q}}\geq \dfrac{2}{1+\dfrac{q_{n-1}}{q}}.
					\end{align*}
					So to prove \eqref{eq.t6.case2a}, it is enough to prove
					$$\dfrac{2}{1+\dfrac{q_{n-1}}{q}}>\dfrac{1}{1-\dfrac{1}{2q}}$$
					which can be written as
					$$1>\dfrac{1+q_{n-1}}{q}.$$
					This is obviously true since $$q=(cb_n+d)q_{n-1}+cq_{n-2}\geq 3q_{n-1}+2q_{n-2}>1+q_{n-1}.$$
					
					\item Suppose that $a_n>b_n,c=1$ and $d=0$. Hence $b_{n+1}$ does not exist and $b_n\geq2$. Then $q=b_nq_{n-1}+q_{n-2}$. Therefore \eqref{eq.t6.case2a} has the form
					\begin{equation} \label{eq.t6.case2a.ii}
						\dfrac{1}{\dfrac{1}{a_n-b_n+r}+\dfrac{q_{n-1}}{q}}>\dfrac{1}{1-\dfrac{1}{2q}}.
					\end{equation}
					\begin{enumerate}
						\item Let $a_n \geq b_n+4$. From $q=b_nq_{n-1}+q_{n-2}\geq b_n \geq 2$ we obtain that $\dfrac{1}{1-\dfrac{1}{2q}} \leq \dfrac{4}{3}$ and we have
						\begin{align*}
							\dfrac{1}{\dfrac{1}{a_n-b_n+r}+\dfrac{q_{n-1}}{q}} &\geq \dfrac{1}{\dfrac{1}{4}+\dfrac{1}{b_n+\dfrac{q_{n-2}}{q_{n-1}}}}\\
							&\geq \dfrac{1}{\dfrac{1}{4}+\dfrac{1}{2}}=\dfrac{4}{3}\geq \dfrac{1}{1-\dfrac{1}{2q}}.
						\end{align*}
						The equality occurs in the third exception.
						
						\item Assume that $n\geq2$ and $a_n=b_n+3$. Therefore we have
						$$\dfrac{1}{\dfrac{1}{a_n-b_n+r}+\dfrac{q_{n-1}}{q}}=\dfrac{1}{\dfrac{1}{3+r}+\dfrac{q_{n-1}}{q}}\geq \dfrac{1}{\dfrac{1}{3}+\dfrac{q_{n-1}}{q}}.$$
						So to prove \eqref{eq.t6.case2a.ii}, it is enough to prove
						$$\dfrac{1}{\dfrac{1}{3}+\dfrac{q_{n-1}}{q}} > \dfrac{1}{1-\dfrac{1}{2q}}$$
						which is equivalent to
						$$\dfrac{4}{3}>\dfrac{1+2q_{n-1}}{q}.$$
						This inequality is true obviously since
						$$q=b_nq_{n-1}+q_{n-2}\geq2q_{n-1}+1.$$
						
						\item Suppose that $n=1$ and $a_1=b_1+3$. Then we have $q=b_nq_{n-1}+q_{n-2}=b_1$. From this we obtain that
						$$\dfrac{1}{\dfrac{1}{3+r}+\dfrac{q_{n-1}}{q}}=\dfrac{1}{\dfrac{1}{3+r}+\dfrac{1}{b_1}}.$$
						So to prove \eqref{eq.t6.case2a.ii} it is enough to prove
						$$\dfrac{1}{\dfrac{1}{3+r}+\dfrac{1}{b_1}}>\dfrac{1}{1-\dfrac{1}{2q}}.$$
						$\bullet$ Let $b_1\geq3$. Then $q\geq3$ and $\dfrac{1}{1-\dfrac{1}{2q}} \leq \dfrac{6}{5}$. Hence
						$$\dfrac{1}{\dfrac{1}{3+r}+\dfrac{1}{b_1}}\geq\dfrac{1}{\dfrac{1}{3}+\dfrac{1}{3}}=\dfrac{3}{2}>\dfrac{6}{5}\geq \dfrac{1}{1-\dfrac{1}{2q}}.$$
						$\bullet$ Suppose that $b_1=2$. Then $q=2$ and $\dfrac{1}{1-\dfrac{1}{2q}} = \dfrac{4}{3}$. It yields
						$$\dfrac{1}{\dfrac{1}{3+r}+\dfrac{1}{b_1}}=\dfrac{1}{\dfrac{1}{3+r}+\dfrac{1}{2}}.$$
						$\bigstar$ Let $a_2$ does not exist then $r=0$. Hence
						$$\dfrac{1}{\dfrac{1}{3+r}+\dfrac{1}{2}}=\dfrac{1}{\dfrac{1}{3}+\dfrac{1}{2}}=\dfrac{6}{5}<\dfrac{4}{3}=\dfrac{1}{1-\dfrac{1}{2q}}.$$
						In this case we have $\alpha=[a_0;5]$ and $\dfrac{p}{q}=[a_0;2]$. This is the fourth exception when $a_2$ does not exist.\\
						$\bigstar$ Suppose that $a_2\geq1$. Then we have
						$$\dfrac{1}{\dfrac{1}{3+r}+\dfrac{1}{2}}<\dfrac{1}{\dfrac{1}{3+1}+\dfrac{1}{2}}=\dfrac{4}{3}=\dfrac{1}{1-\dfrac{1}{2q}}.$$
						This is the fourth exception.
						
						\item Let $n\geq2$ and $a_n=b_n+2$. From this we obtain that
						$$\dfrac{1}{\dfrac{1}{a_n-b_n+r}+\dfrac{q_{n-1}}{q}}\geq\dfrac{1}{\dfrac{1}{2}+\dfrac{q_{n-1}}{q}}.$$
						So to prove \eqref{eq.t6.case2a.ii} it is enough to prove that
						$$\dfrac{1}{\dfrac{1}{2}+\dfrac{q_{n-1}}{q}}\geq \dfrac{1}{1-\dfrac{1}{2q}}$$
						which can be written as
						$$1\geq \dfrac{1+2q_{n-1}}{q}.$$
						This is obviously true since
						$q=b_nq_{n-1}+q_{n-2}\geq 2q_{n-1}+1.$\\
						The equality occurs in the fifth exception.
						
						\item Assume that $n=1$ and $a_n=b_n+2$. \\
						Hence $q=b_nq_{n-1}+q_{n-2}=b_1$. From this we obtain that
						$$\dfrac{1}{\dfrac{1}{a_n-b_n+r}+\dfrac{q_{n-1}}{q}}=\dfrac{1}{\dfrac{1}{2+r}+\dfrac{1}{b_1}}.$$
						$\bullet$ Let $b_1\geq3$. Then $q=b_1\geq3$ and $\dfrac{1}{1-\dfrac{1}{2q}} \leq \dfrac{6}{5}$. It yields
						$$\dfrac{1}{\dfrac{1}{2+r}+\dfrac{1}{b_1}}\geq \dfrac{1}{\dfrac{1}{2}+\dfrac{1}{3}}=\dfrac{6}{5}=\dfrac{1}{1-\dfrac{1}{2q}}.$$
						The equality occurs in the sixth exception.
						
						$\bullet$ Suppose that $b_1=2$. Then $q=b_1=2$ and $\dfrac{1}{1-\dfrac{1}{2q}} = \dfrac{4}{3}$. Hence
						$$\dfrac{1}{\dfrac{1}{2+r}+\dfrac{1}{b_1}}=\dfrac{1}{\dfrac{1}{2+r}+\dfrac{1}{2}}.$$
						$\bigstar$ Assume that $a_2$ does not exist. Then $r=0$ and we obtain that
						$$\dfrac{1}{\dfrac{1}{2+r}+\dfrac{1}{2}}=\dfrac{1}{\dfrac{1}{2}+\dfrac{1}{2}}=1<\dfrac{4}{3}=\dfrac{1}{1-\dfrac{1}{2q}}.$$
						This is the seventh exception when $a_2$ does not exist.\\
						$\bigstar$ Let $a_2\geq1$. It implies that
						$$\dfrac{1}{\dfrac{1}{2+r}+\dfrac{1}{2}}<\dfrac{1}{\dfrac{1}{2+1}+\dfrac{1}{2}}=\dfrac{6}{5}<\dfrac{4}{3}=\dfrac{1}{1-\dfrac{1}{2q}}.$$
						This is the seventh exception.
					\end{enumerate}
				\end{enumerate}
				\item Suppose that $b_n> a_n$ then \eqref{eq.t6.case2} has the form
				\begin{equation} \label{eq.t6.case2b}
					\dfrac{c}{\dfrac{1}{c(b_n-a_n-r)+d}-\dfrac{q_{n-1}}{q}}>\dfrac{1}{1-\dfrac{1}{2q}}.
				\end{equation}
				\begin{enumerate}
					\item Suppose that $c\geq2, d\geq1$ and $b_n\geq a_n+1$. It yields
					\begin{align*}
						\dfrac{c}{\dfrac{1}{c(b_n-a_n-r)+d}-\dfrac{q_{n-1}}{q}}>\dfrac{c}{\dfrac{1}{d}-\dfrac{q_{n-1}}{q}}\geq \dfrac{2}{1-\dfrac{1}{q}}>\dfrac{1}{1-\dfrac{1}{2q}}. 
					\end{align*}
					which is obviously true and inequality \eqref{eq.t6.case2b} follows.
					
					\item Let $c=1,d=0$ and $b_n \geq a_n+2$. From this we obtain that
					\begin{align*}
						\dfrac{1}{\dfrac{1}{b_n-a_n-r}-\dfrac{q_{n-1}}{q}}>\dfrac{1}{\dfrac{1}{2-1}-\dfrac{q_{n-1}}{q}}\geq \dfrac{1}{1-\dfrac{1}{q}}>\dfrac{1}{1-\dfrac{1}{2q}}.
					\end{align*}
					which is obviously true and inequality \eqref{eq.t6.case2b} follows.
				\end{enumerate}
				
			\end{enumerate}
		\end{enumerate}
		The proof of Theorem \ref{hanthot6} is complete.
	\end{proof}
	\begin{proof}[Proof of Example \ref{example1}]
		We have
		$$\left|\alpha-\dfrac{p}{q} \right|=\displaystyle\sum_{n=1}^{\infty}\dfrac{1}{2^{2^n-1}A^{2^n}}-\displaystyle\sum_{n=1}^{N}\dfrac{1}{2^{2^n-1}A^{2^n}}=\displaystyle\sum_{n=N+1}^{\infty}\dfrac{1}{2^{2^n-1}A^{2^n}}.$$
		At the same time
		\begin{align*}
			\dfrac{1}{2q^2}&=\dfrac{1}{2\left(2^{2^N-1}A^{2^N}\right)^2}=\dfrac{1}{2.2^{2^{N+1}-2}A^{2^{N+1}}}\\
			&=\dfrac{1}{2^{2^{N+1}-1}A^{2^{N+1}}}<\displaystyle\sum_{n=N+1}^{\infty}\dfrac{1}{2^{2^n-1}A^{2^n}}=\left|\alpha-\dfrac{p}{q} \right|.
		\end{align*}
		Hence we cannot use Legendre's theorem. On the other side we have
		\begin{align*}
			\dfrac{1}{q^2\left(2-\dfrac{1}{q}\right)}&=\dfrac{1}{\left(2^{2^N-1}A^{2^N}\right)^2\left(2-\dfrac{1}{2^{2^N-1}A^{2^N}}\right)}=\dfrac{1}{\left(2^{2^N-1}A^{2^N}\right)\left(2^{2^N}A^{2^N}-1\right)}\\
			&=\dfrac{1}{2^{2^{N+1}-1}A^{2^{N+1}}}.\displaystyle\sum_{n=0}^{\infty}\dfrac{1}{2^{n2^N}A^{n2^N}}=\displaystyle\sum_{n=0}^{\infty}\dfrac{1}{2^{(n+2)2^N-1}.A^{(n+2)2^N}}\\
			&>\displaystyle\sum_{n=N+1}^{\infty}\dfrac{1}{2^{2^n-1}A^{2^n}}=\left|\alpha-\dfrac{p}{q} \right|.
		\end{align*}
		Thus from Theorem \ref{hanthot2} we obtain that $\dfrac{p}{q}$ is a convergent of $\alpha$.
	\end{proof}
	\begin{proof} [Proof of Example \ref{example2}]
		Example \ref{example1} is an immediate consequence of Example \ref{example2} when we set $A=1$.
	\end{proof}
	\begin{proof} [Proof of Example \ref{example4}]
		We have
		$$\left|\alpha-\dfrac{p}{q} \right|=\displaystyle\sum_{n=1}^{\infty}\dfrac{1}{2^{2^n}A^{2^n}}-\displaystyle\sum_{n=1}^{N}\dfrac{1}{2^{2^n}A^{2^n}}=\displaystyle\sum_{n=N+1}^{\infty}\dfrac{1}{2^{2^n}A^{2^n}}.$$
		At the same time
		$$\dfrac{1}{q^2}=\dfrac{1}{2^{2.2^N}A^{2.2^N}}<\left|\alpha-\dfrac{p}{q} \right|$$
		Hence we cannot use Barbolosi and Jager's theorem. On the other side we have
		\begin{multline*}
			\dfrac{1}{\left(1-\dfrac{1}{2q}\right)q^2}=\dfrac{1}{2^{2.2^N}A^{2.2^N}}\displaystyle\sum_{n=0}^{\infty}\dfrac{1}{\left(2.2^{2^N}A^{2^N}\right)^n}\\
			=\displaystyle\sum_{n=0}^{\infty}\dfrac{1}{\left(2^{(n+2)2^N+n}A^{(n+2)2^N}\right)^n}>\left|\alpha-\dfrac{p}{q} \right|.
		\end{multline*}
		Therefore, from Theorem \ref{hanthot6} we obtain that $\dfrac{p}{q}$ is a convergent or nearest mediant of $\alpha$.
		
	\end{proof}

	\section{Data Availibility}
	
	Data sharing is not applicable to this article as no new data were created or analyzed in this study.

	\section{Declaration}
	
	The authors declare that they have no conflict of interest. 
	
	\section{Acknowledgement}
	Tho Phuoc Nguyen is supported by grant SGS01/P\v{r}F/2024.

	AMS Class: 11J82, 11A55.\\
	Key words and phrases: continued fraction, approximation, Theorem of Legendre.\\
	Jaroslav Han\v{c}l, Tho Phuoc Nguyen, Department of Mathematics, Faculty of Sciences, University of Ostrava, 30.~dubna~22, 701~03 Ostrava~1, Czech Republic.\\
	e-mail: jaroslav.hancl@seznam.cz, phuocthospt@gmail.com\\


\begin{thebibliography}{99}
		
		
		
		\bibitem{barbolosi1988} D. Barbolosi : Fractions continues \`a quotients partiels impairs, Th\`ese, Universit\'e de Provence, Marseille (1988).
		\bibitem{barbolosi1994} D. Barbolosi, H. Jager : On a Theorem of Legendre in the theory of continued fractions, Journal de Th\'eories des Nombres de Bordeaux, Vol. 6, No. 1 (1994), 81--94.
		\bibitem{bh} S. Bahnerov\' a, J. Han\v cl : Sharpening of the theorem of Vahlen and related theorems, J. Ramanujan. Math. Soc., vol. 36, no. 2, (2021), 109--121. 
		\bibitem{billingsley} P. Billingsley : Ergodic Theory and Information, John Wiley and Sons, New York, London, Sydney (1965). 
		\bibitem{borel2} \' E. Borel: Contribution \`a l'analyse arithm\' etique du continu, J. Math. Pures 9, vol. 5, (1903), 329--375.
		\bibitem{borwein} J. Borwein, P. Borwein, Pi and the AGM: A Study in Analytic Number Theory and Computational Complexity, John Wiley \& Sons, 
		New York, 1987.
		\bibitem{Dirichlet} L. G. P. Dirichlet (1842): Verallgemeinerung eines Satzes aus der Lehre von
		den Kettenbruchen nebst einige Anwendungen auf die Theorie der Zahlen. S.-B. Preuss. Akad. Wiss., 93–95. 
		\bibitem{feldman} N. I. Fel'dman, Yu. V. Nesterenko, Transcendental Numbers, Encyclopaedia of Mathematical Sciences, 
		vol. 44: Number Theory IV, A. N. Parshin and I. R. Shafarevich, eds., Springer-Verlag, New York, 1998. 
		\bibitem{grace} J. H. Grace : The classification of rational approximations, Proc. London Math. Soc. 17 (1918), 247--258.
		\bibitem{hancl1} J. Han\v cl, Sharpening of theorems of Vahlen and  Hurwitz and approximation properties of the golden ratio, Arch. Math. (Basel) 105, no. 2, (2015), 129--137.  
		\bibitem{hancl2} J. Han\v cl, Second basic theorem of Hurwitz, Lithuanian Mathematical Journal, vol. 56, no. 1, (2016), 72--76.  
		\bibitem{hancl3} J. Han\v cl, On a Theorem of A. A. Markoff, Results in Math., vol. 76, no. 4, (2021), Article 192. 
		\bibitem{hardy} G. H. Hardy,  E. M. Wright, An introduction to the theory of numbers. Sixth edition. Revised by D. R. Heath-Brown and 
		J. H. Silverman. With a foreword by Andrew Wiles. Oxford University Press, Oxford, 2008. 
		\bibitem{hensley} D. Hensley, Continued fractions, Word Scientific Publishing, (2006). 
		\bibitem{hurwitz} A. Hurwitz, \"Uber die angen\" aherte Darstellung der Irrationalzahlen durch rationale Br\"uche, (German) Math. Ann. 39, 
		no. 2, (1891), 279--284. 
		\bibitem{Ito1988} S. Ito : On Legendre's theorem related to Diophantine approximations, S\'eminaire de Th\'eorie des Nombres, Bordeaux, expos\'e 44 (1987-1988),44-01-44-19.
		\bibitem{Itoosaka} S. Ito : Algorithms with mediant convergents and their metrical theory, Osaka J. Math, (1989), 557-578.
		\bibitem{Ito1985} S. Ito, H. Nakada : On natural extensions of transformations related to diophantine approximations, Number theory and Combinatorics, Word Scientific Pub. (1985), 185-207.
		\bibitem{jager&kraaikamp} H. Jager, C. Kraaikamp : On the approximation by continued fractions, Indag. Math. 51 (1989),289--307.
		\bibitem{jones} W. B. Jones, W. J. Thron, Continued fractions analytic theory and applications, Cambridge University 
		Press, Encyclopedia of Mathematics and its applications 11, (1984). 
		\bibitem{karpenkov} O. Karpenkov, \textit{Geometry of continued fractions.} Algorithms and Computation in Mathematics, 26. Springer, Heidelberg 2013. 
		\bibitem{khinchin} A. Ya. Khinchin, Continued fractions, The University of Chicago Press, Chicago, (1964).
		\bibitem{koksmatheorem} J. F. Koksma : Bewijs van een stelling over kettingbreuken, Mathematica A 6 (1937), 226--233.
		\bibitem{koksma1951} J. F. Koksma : On continued fractions, Simon Stevin 29 (1951/52), 96--102.
		\bibitem{koksma1936} J. F. Koksma : Diophantische Approximationen, Julius Springer, Berlin (1936).
		\bibitem{kraaikamp} C. Kraaikamp : A new class of continued fractions, Acta Arith. 57 (1991), 1--39. 
		\bibitem{legendre} A. M. Legendre (1830). \textit{Th\'eorie des Nombres}, troisi\`eme \'edition, Tome 1. Paris. 
		\bibitem{nakada} H. Nakada : Metrical theory for a class of continued fractions transformations and their natural extensions, Tokyo J. Math. 4 (1981), 399--426. 
		\bibitem{rosen} K. H. Rosen, Elementary number theory and its applications, Addison Wesley, fifth edition, (2005). 
		\bibitem{schmidt} W. Schmidt, Diophantine approximation, Lecture Notes in Mathematics 785, Springer, Berlin, (1980). 
		\bibitem{vahlen} K. Th. Vahlen, \" Uber N\" aherungswerte und Kettenbr\" uche, J. Reine Angew. Math. 115, (1895), 221--233.  
		\bibitem{wall} H. S. Wall, Analytic theory of continued fractions, New York, Chelsea, (1948). 
	\end{thebibliography}
\end{document}